\documentclass[a4paper,10pt, final]{amsart}
\usepackage[english]{certus}
\usepackage{flexisym,breqn}
\usepackage{comment}
\usepackage{amsmath, amssymb, amsthm}
 \usepackage[notref,notcite,draft]{showkeys} 
 \usepackage{bm}

 \newcommand{\op}{{\operatorname{op}}}
 \newcommand{\htens}{\overline{\otimes}}
\newcommand{\CC}{\mathbb{C}}
 \newcommand{\ZZ}{\mathbb{Z}}
 \newcommand{\NN}{\mathbb{N}}
 \newcommand{\FF}{\mathbb{F}}
\newcommand{\MM}{\mathbb{M}}
\newcommand{\T}{{\operatorname{(T)}}}
\newcommand{\Der}{{\operatorname{Der}}}
\newcommand{\Pol}{{\operatorname{Pol}}}
\newcommand{\tens}{\otimes}

\newcommand{\Tor}{{\operatorname{Tor}}}
\newcommand{\Ext}{{\operatorname{Ext}}}
\newcommand{\bet}{\beta^{(2)}}
\newcommand{\twoone}{\operatorname{II}_1}

\theoremstyle{definition}
\newtheorem{Por}[Def]{Porism}

\newcommand{\Proj}{{\operatorname{Proj}}}
\newcommand{\ind}{{\operatorname{ind}}}
\newcommand{\sa}{{\operatorname{sa}}}

\usepackage{color}

\title[Measure continuous derivations]{\texorpdfstring{Measure continuous derivations on von Neumann algebras and applications to $\bm{L^2}$-cohomology}{Measure continuous derivations on von Neumann algebras and applications to {L2}-cohomology}
}

\author{Vadim Alekseev} 
\address{Vadim Alekseev,
Mathematisches Institut,
Georg-Au\-gust-Uni\-versi\-t{\"a}t G{\"o}t\-ting\-en,
Bunsenstra{\ss}e 3-5,
D-37073 G{\"o}ttingen, 
Germany.}
\email{alekseev@uni-math.gwdg.de}

\author{David Kyed} 
\address{David Kyed,
Department of Mathematics,
KU Leuven,
Celestijnenlaan 200B,
B-3001 Leuven, 
Belgium.}
\email{David.Kyed@wis.kuleuven.be}
\urladdr{www.kuleuven.be/~u0078326}

\keywords{von Neumann algebras, $L^2$-Betti numbers,  property $\T$. }
\subjclass[2010]{46L10, 46L52, 46L57}

\thanks{The second named author gratefully acknowledges the funding from The Danish Council for Independent Research $|$ Natural Sciences and the ERC Starting Grant VNALG-200749}

\begin{document}

\begin{abstract}
We prove that norm continuous derivations from a von Neumann algebra into the algebra of operators affiliated with its tensor square are automatically continuous for both the strong operator topology and the measure topology. Furthermore, we prove that the first continuous $L^2$-Betti number scales quadratically when passing to corner algebras and  derive an upper bound given by Shen's generator invariant. This, in turn, yields vanishing of the first continuous $L^2$-Betti number for $\twoone$ factors with property (T), for finitely generated factors with non-trivial fundamental group and for factors with property Gamma.

\end{abstract}

\maketitle

\section*{Todo}

%

\section{Introduction}
The theory of $L^2$-Betti numbers has been generalized to a vast number of different contexts since the seminal work of Atiyah \cite{atiyah-l2}. One  recent such generalization is due to Connes and Shlyakhtenko \cite{CS} who introduced $L^2$-Betti numbers for subalgebras of finite von Neumann algebras, with the main purpose being to obtain a suitable notion  for arbitrary II$_1$-factors. Although their definitions are very natural, it has proven to be quite difficult to perform concrete calculations. The most advanced computational result so far is due to Thom \cite{thom2008l2} who proved that the $L^2$-Betti numbers vanish for von Neumann algebras with diffuse center. Notably, the problem of computing a positive degree $L^2$-Betti number for a single II$_1$-factor has remained open for a decade at the time of writing! Due to this evident drawback, Thom \cite{thom2008l2} introduced a continuous version of the first $L^2$-Betti number, which 
turns out to be much more manageable than its algebraic counterpart. The first continuous $L^2$-Betti number is defined as the von Neumann dimension of the first continuous Hochschild cohomology of the von Neumann algebra $M$ with values in the algebra of operators affiliated with $M\htens M^\op$. The word `continuous' here means that we restrict attention to those derivations that are continuous from the norm topology on $M$ to the measure topology on the affiliated operators.\\

In this paper we continue the study of Thom's continuous version of the first $L^2$-Betti number and our first result (Theorem \ref{sot-mt-cont}) shows that norm continuous derivations are automatically also continuous for both the strong operator topology and the measure topology. This allows us  to derive all previously known computational results concerning the first continuous $L^2$-Betti number, and furthermore to prove that it vanishes for $\twoone$ factors with property $\T$ (Theorem \ref{beta-1-of-T-factor}). In Section \ref{compression-subsection} we prove that it scales quadratically when passing to corner algebras (Theorem \ref{scaling-for-eta}) and is dominated by Shen's generator invariant (Corollary \ref{bounded-by-G-cor}). Along the way, we give a new short cohomological proof of the fact that the (non-continuous) first $L^2$-Betti number vanishes for von Neumann algebras with diffuse center,  and furthermore derive a number of new vanishing results regarding the first continuous $L^2$-Betti 
number, including the vanishing for $\twoone$ factors with property Gamma and finitely generated factors with non-trivial fundamental group (Corollaries \ref{bounded-by-G-cor} \& \ref{fundamenal-grp-cor}).

\section{Preliminaries}
In this section we briefly recapitulate the theory of non-commutative integration and the theory of $L^2$-Betti numbers for von Neumann algebras.

\subsection{Non-commutative integration} 
Let us recall some facts from the theory of non-commutative integration, cf. \cite{nelson-nc-integration}, \cite[IX.2]{takesaki-2}.
Let $N$ be a finite von Neumann algebra equipped with a normal, faithful, tracial state $\tau$. Consider $N$ in its representation on the GNS-space arising from $\tau$, and let $\ms N$ be the algebra of (potentially) unbounded, closed,  densely defined operators affiliated with $N$. We equip $\ms N$ with the \emph{measure topology}, defined by the following two-parameter family of neighbourhoods of zero:
\beqn
N(\eps,\delta) = \{a\in \ms N\,|\, \exists p\in \mathrm{Proj}(N): \norm{ap} < \eps,\; \tau(p^\perp) < \delta\}, \quad \eps,\delta>0.
\eeqn
With this topology, $\ms N$ is a complete \cite[Theorem IX.2.5]{takesaki-2}  metrizable \cite[Theorem 1.24]{rudin-funct-an}  topological vector space  and the multiplication map
\beqn
(a,b)\mapsto ab \colon \ms N\times \ms N \to \ms N 
\eeqn
is uniformly continuous when restricted to products of bounded subsets \cite[Theorem 1]{nelson-nc-integration}. Convergence with respect to the measure topology is also referred to as \emph{convergence in measure}.
We also introduce the notation
\beqn
N(0,\delta) = \{a\in \ms N\,|\, \exists p\in \mathrm{Proj}(N): ap = 0,\; \tau(p^\perp) < \delta\},
\eeqn
and
\beqn
N(\eps, 0) = \{a\in N\,|\, \norm{a} < \eps \}\subset \ms N.
\eeqn
Notice that $N(0,\delta)$ and $N(\eps,0)$ are \emph{not} zero neighbourhoods in the measure topology, but merely $G_\delta$ sets
However, the following additive and multiplicative properties continue to hold for all $\eps_1,\eps_2,\delta_1,\delta_2 \geqslant 0$, cf. \cite[Theorem 1]{nelson-nc-integration}:
\beq\label{eq:measure-top-addition}
N(\eps_1,\delta_1) + N(\eps_2,\delta_2) \subset N(\eps_1 + \eps_2,\delta_1+\delta_2),
\eeq
\beq\label{eq:measure-top-multiplication}
N(\eps_1,\delta_1) \cdot N(\eps_2,\delta_2) \subset N(\eps_1\eps_2,\delta_1+\delta_2).
\eeq
The noncommutative $L^p$-spaces $L^p(N,\tau)$ are naturally identified with subspaces of $\ms N$ \cite[Theorem IX.2.13]{takesaki-2}. We fix the notation $\xrightarrow{s}$ for strong convergence of elements in von Neumann algebras, $\xrightarrow{2}$ for the $L^2$-convergence and $\xrightarrow{m}$ for the convergence in measure of elements in $\ms N$. Clearly strong convergence implies convergence in 2-norm, and we remind the reader that for nets that are bounded in the operator norm the converse is also true --- a fact we will use extensively in the sequel. As in the commutative case, the Chebyshev inequality can be used to establish the following fact.
\begin{Lemma}[{\cite[Theorem 5]{nelson-nc-integration}}]\label{Lp-lem}
For any $p\geqslant 1$ the inclusion $L^p(N,\tau)\subset \ms N$ is continuous; i.e.~$L^p$-convergence implies convergence in measure.
\end{Lemma}

\subsection{\texorpdfstring{$\bm{L}^{\bm2}$-Betti numbers for tracial algebras}{L2-Betti numbers for tracial algebras}}
In \cite{CS} Connes and Shlyakh\-tenko introduced $L^2$-Betti numbers in the general setting of tracial $*$-algebras; if $M$ is a finite von Neumann algebra and $\mc A\subset M$ is any weakly dense unital $*$-subalgebra its $L^2$-Betti numbers are defined as
\[
\beta_p^{(2)}(\mc A, \tau)=\dim_{M\htens M^\op} \Tor_p^{\mc A \odot \mc A^\op }(M\htens M^\op, \mc A).
\]
Here the dimension function $\dim_{M\htens M^\op}(-)$ is the extended von Neumann dimension due to L{\"u}ck; cf.~\cite[Chapter 6]{Luck02}.
This definition is inspired by the well-known correspondence between representations of groups and bimodules over finite von Neumann algebras, and it extends the classical theory by means of the formula $\beta_p^{(2)}(\Gamma)=\beta_p^{(2)}(\mathbb{C}\Gamma,\tau)$ whenever $\Gamma$ is a discrete countable group. In \cite{thom2008l2} it is shown that the $L^2$-Betti numbers also allow the following cohomological description:
\[
\beta_p^{(2)}(\mc A,\tau)=\dim_{M\htens M^\op}\Ext_{\mc A \odot \mc A}^p(\mc A, \ms U),
\]
where $\ms U$ denotes the algebra of operators affiliated with $M\htens M^\op$. It is a classical fact \cite[1.5.8]{loday} that the Ext-groups above are isomorphic to the Hochschild cohomology groups of $\mc A$ with coefficients in  $\ms U$, where the latter is considered as an $\mc A$-bimodule with respect to the actions
\[
a\cdot \xi:= (a\tens1^\op)\xi \ \text{ and } \ \xi\cdot b:= (1\tens b^\op)\xi \ \text{ for } a,b \in \mc A \text{ and } \xi \in \ms U. 
\]
In particular, the first $L^2$-Betti number can be computed as the dimension of the right $M\htens M^\op$-module
\[
H^1(\mc A, \ms U)=\frac{\Der(\mc A, \ms U)}{\text{Inn}(\mc A,\ms U)}. 
\]
Here $\Der(\mc A,\ms U) $ denotes the space of derivations from $\mc A$ to $\ms U$ and $\text{Inn}(A,\ms U)$ denotes the space of inner derivations. We recall that a linear map $\delta $ from $\mc A$ into an $\mc A$-bimodule $\mc X$ is called a \emph{derivation} if it satisfies
\[
\delta(ab)=a\cdot \delta(b) + \delta(a)\cdot b \ \text{ for all } a,b\in \mc A,
\]
and that a derivation is called \emph{inner} if there exists a vector $\xi \in \mc X$ such that
\[
\delta(a)=a\cdot \xi-\xi\cdot a \text{ for all } a\in \mc A.
\]
When the bimodule in question is $\ms U$, with the bimodule structure defined above, the derivation property amounts to the following:
\[
\delta(ab)=(a\otimes 1^\op)\delta(b) + (1\otimes b^\op)\delta(a) \ \text{ for all } a,b\in \mc A.
\]
Although the extended von Neumann dimension is generally not faithful, enlarging the coefficients from $M\htens M^\op$ to $\ms U$  has the effect that $\beta_1^{(2)}(\mc A,\tau)=0$ if and only if $H^1(\mc A, \ms U)$ vanishes \cite[Corollary 3.3 and Theorem 3.5]{thom2008l2}. In particular, in order to prove that $\beta_1^{(2)}(\mc A, \tau)=0$ one has to prove that every derivation from $\mc A$ into $\ms U$ is inner.  
These purely algebraically defined $L^2$-Betti numbers have turned out extremely difficult to compute in the case when $\mc A$ is $M$ itself. Actually, the only computational result  known in this direction (disregarding finite dimensional algebras) is that the they vanish for von Neumann algebras with diffuse center  (see \cite[Corollary 3.5]{CS} and \cite[Theorem 2.2]{thom2008l2}). In particular, for $\twoone$-factors not a single computation of a positive degree $L^2$-Betti is known, and furthermore this seems out of reach with tools available at the moment. It is therefore natural to consider variations of the definitions above that take into account the topological nature of $M$, and in \cite{thom2008l2} Thom suggests to consider a first cohomology group consisting of (equivalence classes of) those derivations $\delta\colon \mc A \to \ms U$ that are closable from the norm topology to the measure topology. Note that when $\mc A$ is norm closed these are exactly the derivations that are norm-measure 
topology continuous. We denote the space of closable derivations by $\Der_c(\mc A, \ms U)$, the  continuous cohomology by $H^1_c(\mc A, \ms U)$ and by $\eta_1^{(2)}(\mc A,\tau)$ the corresponding continuous $L^2$-Betti numbers; i.e.~ 
\[
\eta_1^{(2)}(\mc A,\tau)=\dim_{M\htens M^\op} H^1_c(\mc A,\ms U).
\]
These continuous $L^2$-Betti numbers are much more manageable than their algebraic counterparts --- for instance they are known \cite[Theorem 6.4]{thom2008l2} to vanish for von Neumann algebras that are non-prime and for those that contain a diffuse Cartan subalgebra.  \\

Finally, let us fix a bit of notation. For the rest of this paper, we consider a finite von Neumann algebra $M$ with separable predual $M_*$. We endow $M$ with a fixed faithful, normal, tracial state $\tau$ and consider $M$ in the GNS representation on the Hilbert space $\mc H = L^2(M,\tau)$. The trace $\tau$ induces a faithful, normal, tracial state on the von Neumann algebraic tensor product $M\overline{\otimes} M^{\textup{op}}$ of $M$ with its opposite algebra; abusing notation slightly, we will still denote it by $\tau$. We always consider  $M\overline{\otimes} M^{\textup{op}}$ in the GNS representation  on $L^2(M\overline{\otimes}M^{\textup{op}},\tau)$ and denote by $\ms U$ the algebra of closed, densely defined, unbounded operators affiliated with $M\overline{\otimes} M^{\textup{op}}$. More generally, if $N$ is a finite von Neumann algebra endowed with a tracial state $\rho$ we denote by $\ms U (N)$ the algebra of operators on $L^2(N,\rho)$ affiliated with $N$. We will use the symbol ``$\htens$'' to 
denote tensor products of von Neumann algebras and ``$\odot$'' to denote algebraic tensor products, and, unless explicitly stated otherwise, all subalgebras in $M$ are implicitly assumed to contain the unit of $M$. Finally, when there is no source of confusion we will often suppress the notational reference to the trace $\tau$, and simply write $L^2(M), \beta_p^{(2)}(M), \eta_1^{(2)}(M)$ etc.

\section{Improving the continuity}

In this section we prove that a derivation which is continuous for the norm topology is automatically continuous for the strong operator topology and the measure topology as well. Intuitively, this statement is based on the fact that strong convergence is ``almost uniform'', which is known as the non-commutative Egorov theorem, cf.~\cite[Theorem II.4.13]{takesaki-1}. The precise statement is as follows.

\begin{Thm}\quad\label{sot-mt-cont}
\begin{enumerate}
 \item Let $A\subset M$ be a weakly dense \Cs-algebra and let $\delta\colon A\to \ms U$ be a norm-measure continuous derivation. Then $\delta$ has a unique norm-measure continuous extension to $M$.
 \item Let $\delta\colon M\to \ms U$ be a norm-measure continuous derivation. Then $\delta$ is also continuous from the strong operator topology to the measure topology.
 \item Let $\delta\colon M\to \ms U$ be a norm-measure continuous derivation. Then $\delta$ is also continuous from the measure topology on $M$ to the measure to\-po\-lo\-gy on $\ms U$; in particular, it has a unique measure-measure continuous extension to the algebra $\ms M$ of operators affiliated with $M$.
\end{enumerate}
\end{Thm}
Note that the extension property in i) also follows from the rank-metric based arguments in \cite[Theorem 4.3 \& Proposition 4.4]{thom2008l2}, but the boundedness of the extension is not apparent from this.

\begin{proof}
We first prove i). By Kaplansky's density theorem the unit ball $(A)_1$ is strongly dense in $(M)_1$, and for $a\in (M)_1$ we may therefore choose a sequence\footnote{As the predual $M_*$ is assumed separable, $(M)_1$ is a separable and metrizable space for the strong operator topology and we can therefore do with sequences rather than nets.} $a_n\in (A)_1$ with $a_n\xrightarrow{s} a$. Let us first prove that $\delta(a_n)$ is Cauchy in measure. So we fix an $\eps > 0$ and want to find an $N_0$ such that $\delta(a_n)-\delta(a_m) \in N(\eps,\eps)$ for $\min\{m,n\}>N_0$.  We first make use of the norm-measure topology continuity of $\delta\colon A\to \ms U$ to find a $\gamma > 0$ such that $\norm{a}\leqslant \gamma$ implies that $\delta(a)\in N(\eps/3,\eps/3)$. 
Consider now $\mb N \times \mb N$ with the ordering
\[
(m,n)\geqslant (m',n') \text{ iff } m\geqslant m'  \text{ and }  n\geqslant n', 
\]
and the net of self-adjoint elements $b_{(m,n)} := (a_m-a_n)^*(a_m-a_n)$. Then for every $\xi\in \mc H$ we have $\norm{b_{(m,n)}\xi} \leqslant2 \norm{(a_m-a_n)\xi}$ and hence $b_{(m,n)}\xrightarrow{s}0$. 
Let now $f\colon \mb R\to [0,1]$ be a continuous function with $f(x)=1$ for $x\leqslant \gamma^2/4$ and $f(x) = 0$ for $x\geqslant \gamma^2$ and consider the net $h_{(m,n)} := f(b_{(m,n)})$. It follows that $\norm{h_{(m,n)}}\leqslant 1$ and 
\beq
\frac{\gamma^2}{4}(1 - h_{(m,n)}) \leqslant b_{(m,n)}.\notag
\eeq
Since the $C^*$-algebra generated by $b_{(m,n)}$ is commutative and both $b_{(m,n)}$ and $1-h_{(m,n)}$ are positive this implies
\beq
0 \leqslant (1 - h_{(m,n)})^*(1-h_{(m,n)}) \leqslant \tfrac{16}{\gamma^4} \ b_{(m,n)}^*b_{(m,n)}^{\phantom{*}},\notag
\eeq
and thus
\[
\|1-h_{(m,n)}\|_2\leqslant \tfrac{4}{\gamma^2}\|b_{(m,n)}\|_2\rightarrow 0.
\]
Hence $1-h_{(m,n)}\xrightarrow{2} 0$, and the convergence therefore holds in measure as well.
Also note that $\norm{b_{(m,n)}h_{(m,n)}}\leqslant \gamma^2$ and hence

\beq \label{norm-estimate}
\norm{(a_m-a_n)h_{(m,n)}}\leqslant \gamma.
\eeq
We now use the derivation property to obtain:
\begin{dmath*}[breakdepth=0]
\delta(a_m)-\delta(a_n) = \delta((a_m-a_n)h_{(m,n)}) + \delta((a_m-a_n)(1-h_{(m,n)}))
=\delta((a_m-a_n)h_{(m,n)})  +(1\otimes (1-h_{(m,n)})^\op)\delta(a_m-a_n)  \\ - ((a_m-a_n)\otimes1^\op)\delta(h_{(m,n)}).
\end{dmath*}
By \eqref{norm-estimate} and the choice of $\gamma$, the first summand is in $N(\eps/3,\eps/3)$. Let us now consider the second summand. The norm-measure topology boundedness of $\delta$ on $A$ implies \cite[Theorem 1.32]{rudin-funct-an} that the set 
\[
\{\delta(a_m-a_n) \mid n,m\in\mb N \}
\]
 is bounded in $\ms U$. Hence it follows from the fact that $1-h_{(m,n)}\xrightarrow{m} 0$ together with the uniform continuity of multiplication on bounded sets of $\ms U$, that there exists an $N_1$ such that $(1\otimes(1-h_{(m,n)})^\op)\delta(a_m-a_n)\in N(\eps/3,\eps/3)$ for $\min\{m,n\} > N_1$. Lastly we consider the third term.
Again by norm-boundedness of $\delta$, the set $\{\delta(h_{(m,n)})\}_{n,m\in \mb N}$ is bounded in $\ms U$. As $a_n$ is  strongly convergent it also converges in 2-norm; hence $a_n\tens 1^\op$ converges in 2-norm and is, in particular, a Cauchy sequence for the measure topology.  Thus, there exists an $N_2$ such that 
\[
((a_m-a_n)\otimes 1^\op)\delta(h_{(m,n)})\in N(\eps/3,\eps/3) 
\]
for $\min\{n,m\} > N_2$. Taking $N_0 = \max\{N_1,N_2\}$ establishes that $\delta(a_n)$ is a Cauchy sequence in the measure topology. Appealing to the completeness of $\ms U$, we may now define $\delta(a):=\lim_n \delta(a_n)$; it is routine to check that this yields a well-defined derivation $\delta\colon M\to \ms U$. The continuity of the extension follows from its definition from which it is clear that the set $\delta((M)_1)$ is contained in the measure topology closure of the bounded set $\delta((A)_1)$ which is again bounded by \cite[Theorem 1.13]{rudin-funct-an}. Thus the extension maps norm bounded sets to measure topology bounded sets and is therefore continuous \cite[Theorem 1.32]{rudin-funct-an}, and the proof of i) is complete.\\

Next we prove ii) and iii). If $\delta\colon M\to \ms U$ is a norm continuous derivation, then by repeating the above arguments for $m=\infty$, $a=a_\infty$, we get strong continuity of $\delta$ on bounded sets: if $a_n \xrightarrow{s} a$ and the sequence $a_n$ is uniformly bounded, then $\delta(a_n)\xrightarrow{m} \delta(a)$. 
Since strong convergence implies $L^2$-convergence it also implies convergence in measure by Lemma \ref{Lp-lem}, and hence it suffices to prove that $\delta$ is measure-measure continuous. To this end, by metrizability of the measure topology it suffices to take a sequence $a_n \in M$ such that $a_n \xrightarrow{m} 0$ and prove that $\delta(a_n)\xrightarrow{m} 0$. Since $a_n\xrightarrow{m} 0$
 we get\footnote{The existence of $p_n$ can, for instance, be seen by noting that $ \eps_n := \inf \{ \eps>0 \mid  a_n \in N(\epsilon,\epsilon)  \}$ must converge to zero if $a_n\overset{m}{\to}0$. }  a sequence of projections $p_n$ such that $\norm{a_n p_n} \to 0$ and $\tau(p_n)\to 1$
Thus, $p_n$ is a norm bounded sequence that converges to $1$ strong\-ly and, by what was just proven, it follows that $\delta(p_n)\xrightarrow{m} \delta(1)=0$.
Now we use the derivation property of $\delta$:
\beq\label{deriv-lambda-eq}
\delta(a_n) = \delta(a_n p_n) + (1\otimes (1-p_n)^\op)\delta(a_n) - (a_n\otimes 1^\op)\delta(p_n).
\eeq
As $\delta$ is norm-bounded, the first summand in \eqref{deriv-lambda-eq} converges to $0$ in measure. The second summand converges to $0$ in measure because $1-p_n \in N(0,\eps_n)$ with $\eps_n=1-\tau(p_n) \to 0$ and $\delta(a_n)\in N(\gamma_n,\eps_n)$ for some $\gamma_n > 0$ \cite[Lemma IX.2.3]{takesaki-2}. The third summand converges to $0$ in measure because $\delta(p_n)\xrightarrow{m} 0$ and multiplication is measure continuous.  This finishes the proof.
\end{proof}

The next lemma shows that there is no hope for weaker continuity properties of derivations.
\begin{Lemma}
 Let $M$ be a diffuse finite von Neumann algebra. The only derivation $\delta\colon M\to \ms U$ which is continuous from the ultraweak topology on $M$ to the measure topology is the zero map.
\end{Lemma}
\begin{proof}
 Let $\delta\colon M\to \ms U$ be a derivation which is continuous from the ultraweak topology on $M$ to the measure topology. Let $\{u_n\}_{n\in\mb N}\subset M$ be a sequence of unitaries weakly converging to zero and let $m\in M$ be given. Then $u_nm$ weakly converges to zero and we have
 \[
  (u_n\otimes 1^\op)\delta(m)= \delta(u_nm)  - 1\otimes m^\op\delta(u_n).
 \]
Because the ultraweak topology and the weak operator topology agree on bounded subsets of $M$, both summands on the right-hand side converge to zero in measure since  $\delta$ is assumed continuous; hence $(u_n\otimes 1^\op)\delta(m)\xrightarrow{m} 0$. Multiplying by the unitaries $(u_n^*\otimes 1^\op)$, we infer $\delta(m) = 0$.
\end{proof}

\begin{Rem}
Bearing in mind the numerous automatic continuity results for derivations between operator algebras (see \cite{sinclair-smith} and references therein), it is of course natural to ask if norm continuity of a derivation $\delta\colon M\to \ms U$ is also automatic. We were not able to prove this, and it seems to be a difficult question to answer. One reason being the absence of examples of finite von Neumann algebras (or $C^*$-algebras, for that matter) for which a non-inner derivation into the algebra $\ms U$ is known to exist. Moreover, the fact that we are considering the operators affiliated with $M\htens M^\op$ (as opposed to $M$ itself) has to play a role if automatic norm continuity is to be proven, as there are examples of derivations from $M$ into the operators affiliated with $M$ which are not norm-measure topology continuous. This follows from \cite{discontinuous-derivations} where the authors  exhibit a (commutative) finite von Neumann algebra $M$ for which there exists a derivation $\delta\colon \
\ms{U}(M)\to \ms{U}(M)$ which is not measure-measure continuous. If the restriction $\delta|_M $ were norm-measure continuous, then it is not difficult to see that the graph of the original derivation $\delta$ is closed (in the product of the measure topologies), and  hence it cannot be discontinuous. 
If norm continuity turns out to be automatic, there is of course no difference between the ordinary and the continuous $L^2$-Betti numbers, and one might even take the standpoint that if there is no such automatic continuity, then continuity has to be imposed in order to get a satisfactory theory.
\end{Rem}

\section{\texorpdfstring{The first continuous ${L}^{2}$-Betti number}{The first continuous L2-Betti number}}
In this section we apply the above automatic continuity result to obtain information about the first continuous $L^2$-Betti number for von Neumann algebras. Some of the results presented are already known, or implicit in the literature, but since the proofs are knew and quite simple we hope they will shed new light on these results. The main new result in this section is Theorem \ref{beta-1-of-T-factor} which shows that the first continuous $L^2$-Betti number vanishes for property $\T$ factors.\\

Recall that if $N\subset M$ is an inclusion of von Neumann algebras, then the normalizer of $N$ in $M$ is defined as the set of unitaries in $M$ which normalize $N$:
\beqn
\mc N_M(N) = \{u \in U(M)\,|\, u^*Nu = N\}.
\eeqn
The following lemma appears in \cite{thom2008l2}, but we include its short proof for the sake of completeness.
\begin{Lemma}[{\cite[Lemma 6.5]{thom2008l2}}]\label{lemma:vanishing-on-the-normalizer}
 Let $\delta\colon M\to \ms U$ be a derivation which vanishes on a diffuse subalgebra $N\subset M$ and let $u\in \mc N_M(N)$. Then $\delta(u) = 0$.
\end{Lemma}
\begin{proof}
Let $h\in N$ be a diffuse element. Since $\delta(u^*)= - (u^*\tens u^{*\op})\delta(u)$ we get
\begin{align*}
0 = \delta(uhu^*) &= (1\otimes (hu^*)^\op)\delta(u) + (uh\otimes 1^\op)\delta(u^*)\\
&=  (1\otimes(hu^*)^\op)\delta(u) - (uhu^*\otimes u^{*\,\op})\delta(u) \\
&= (u\otimes u^{*\op})(1\otimes h^\op-h\otimes 1^\op)(u^*\otimes 1)\delta(u).
\end{align*}
Since $h$ is diffuse, $1\otimes h^{\textup{op}} - h\otimes 1^\op$ is not a zero divisor in $\ms U$, and it follows that $\delta(u)=0$.

\end{proof}
The next lemma, which might be of independent interest, shows that  if a von Neumann algebra only allows inner derivations into the algebra of operators affiliated with its double, then the same is true for the algebra of operators affiliated with the double of any ambient von Neumann algebra.

\begin{Lemma}\label{lemma:vanishing-change-of-coeff}
Let $N\subset M$ be a sub-von Neumann algebra. Then the following statements are equivalent:
\begin{enumerate}
 \item $\beta_1^{(2)}(N,\tau)=0$,
 \item every derivation $\delta\colon N\to \ms U$ is inner (where $\ms U$ is considered as an $N$-bimodule via the inclusion $N\subset M$).
\end{enumerate}
\end{Lemma}
\begin{proof}
By \cite[1.5.9]{loday} we have $H^1(N, \ms U)=\Ext^1_{N\odot N^\op}(N, \ms U)$  and by \cite[Theorem 3.5]{thom2008l2} this right $\ms U$-module is isomorphic to
\[
\Hom_{\ms U} \left({\Tor^{N\odot N^\op}_1(\ms U, N)},{\ms U}\right).
\]
Furthermore, by \cite[Corollary 3.3]{thom2008l2} this module vanishes exactly when
\[
\dim_{M\htens M^\op} \Tor^{N\odot N^\op}_1(\ms U, N)=0.
\]
But since $\ms U\odot_{M\htens M^\op} -$ and $M\htens M^\op\odot_{N\htens N^\op} -$ are both flat and dimension preserving (see \cite[Proposition 2.1 and Theorem 3.11]{reich-K-and-L-of-aff-op} and \cite[Theorem 6.29]{Luck02}) we get
\begin{align*}
\dim_{M\htens M^\op} \Tor^{N\odot N^\op}_1(\ms U, N) &= \dim_{M\htens M^\op} \Tor_1^{N\odot N^\op}(\ms U\odot_{M\htens M^\op}M\htens M^\op, N)\\
&=\dim_{M\htens M^\op} \ms U\odot_{M\htens M^\op}  \Tor_1^{N\odot N^\op}(M\htens M^\op, N)\\
&=\dim_{M\htens M^\op}  \Tor_1^{N\odot N^\op}(M\htens M^\op, N)\\
&=\dim_{M\htens M^\op}  \Tor_1^{N\odot N^\op}(M\htens M^\op \odot_{N\htens N^\op}N\htens N^\op , N)\\
&=\dim_{M\htens M^\op} M\htens M^\op \odot_{N\htens N^\op}  \Tor_1^{N\odot N^\op}(N\htens N^\op , N)\\
&=\dim_{N\htens N^\op}  \Tor_1^{N\odot N^\op}(N\htens N^\op , N)\\
&=\beta_1^{(2)}(N,\tau). \qedhere
\end{align*} 
\end{proof}

Since $H^1(M,\ms U)$ is the algebraic dual of $H_1(M,\ms U)$ \cite[Theorem 3.5]{thom2008l2} it follows from \cite[Corollary 3.3 \& 3.4]{thom2008l2}  that $\beta_1^{(2)}(M)=0$ if and only if $H^1(M,\ms U)$ is trivial. In the continuous case it is not so clear if the cohomology $H^1_c(M,\ms U)$ is also a dual module, but as the following proposition shows, vanishing of $\eta_1^{(2)}(M)$ does actually imply innerness of all continuous derivations from $M$ to $\ms U$. As we will see in the following subsections, $\eta_1^{(2)}(M)$ vanishes in a lot of special cases and these vanishing results can therefore be translated into automatic innerness of continuous derivations on $M$ with values in $\ms U$. We point out that we do not know an example of a von Neumann algebra $M$ for which $\eta_1^{(2)}(M)\neq 0$ and it cannot be excluded that a innerness of continuous derivations from $M$ to $\ms U$ is automatic in general.

\begin{Prop}
We have $\eta_1^{(2)}(M)=0$ if and only if $H^1_c(M,\ms U)=\{0\}$.
\end{Prop}
The proof is a modification of an argument from \cite{peterson-thom}.
\begin{proof}
The ``if'' part of the statement is obvious so we only have to prove the converse. Assume therefore that $\eta_1^{(2)}(M):=\dim_{M\htens M^\op}H^1_{c}(M,\ms U)=0$ and let $\delta\in \Der_c(M,\ms U)$ be given; we need to prove that $\delta$ is inner. By Sauer's local criterion \cite[Theorem 2.4]{sauer-betti-for-groupoids}  we can find a partition $\{p_n\}_{n=1}^\infty$ of the unit in $M\htens M^\op$ such that $\delta(-)p_n\in \Inn(M,\ms U)$ for every $n\in \NN$. Thus, there exists $\xi_n\in \ms U$ such that
\[
\delta(x)p_n=(x\tens 1-1\tens x^\op)\xi_n \ \text{ for all } x\in M,
\]
and we may therefore furthermore assume that $\xi_n=\xi_np_n$ for every $n\in \NN$. We now claim that the series $\sum_{n=1}^\infty \xi_n$ converges in measure and that its limit implements $\delta$. To prove the convergence it suffices to show that $S_k:=\sum_{i=1}^k \xi_i $ is Cauchy in measure. For given $\eps>0$ and $k$ sufficiently big we have $\sum_{i=k+1}^\infty \tau(p_i)<\eps$, and therefore $q_k:=\sum_{i=1}^k p_i$ satisfies $\tau(1-q)<\eps$ and for every $l\geqslant k$ we have
\[
(S_l-S_k)q= \sum_{i=k+1}^l\sum_{j=1}^k \xi_ip_j =\sum_{i=k+1}^l\sum_{j=1}^k \xi_ip_ip_j =0,
\]
since the $p_n$'s are mutually orthogonal. Thus, $S_k$ is Cauchy in the measure topology (even in the rank topology) and hence the limit $\xi:=\lim_k^m  S_k$ exists and this limit implements $\delta$ since
\begin{align*}
\delta(x)&=\lim^m_{k\to \infty} \delta(x)q_k=\lim_{k\to \infty}^m\sum_{i=1}^k \delta(x)p_i =\lim_{k\to \infty}^m \sum_{i=1}^k (x\tens 1-1\tens x^\op)\xi_i\\
&=\lim_{k\to \infty}^m (x\tens 1-1\tens x^\op)S_k=(x\tens 1-1\tens x^\op)\xi. \qedhere
\end{align*}
\end{proof}
Note that the above proof does not uses the continuity of the derivation $\delta$ at any point, so this also provides a proof of the fact that $\beta_1^{(2)}(M)$ vanishes iff $H^1(M,\ms U)=\{0\}$.

\subsection{The case of diffuse center}
Using methods from free probability, Connes and Shlyakhtenko proved in \cite[Corollary 3.5]{CS} that $\beta_1^{(2)}(M,\tau)=0$  when $M$ has diffuse center, and using homological algebraic methods this was later generalized by Thom to higher $L^2$-Betti numbers in \cite[Theorem 2.2]{thom2008l2}. In this section we give a short cohomological proof of this result in degree one.\\

\begin{Prop}[{\cite[Corollary 3.5]{CS}}]\label{beta-one-of-diffuse}
If $M$ has diffuse center then every deri\-va\-tion $\delta\colon M\to \ms U$ is norm-measure topology continuous and $\beta_1^{(2)}(M,\tau)$ is zero.
\end{Prop}

\begin{proof}
Since the center $Z(M)$ is diffuse we can choose an identification $Z(M)=L^\infty(\mathbb{T})=L\ZZ$. Denote by $h$ a diffuse, selfadjoint generator of $Z(M)$. To see that $\delta$ is bounded, it suffices to prove that its graph is closed \cite[Theorem 2.15]{rudin-funct-an}; let therefore $x_n\in M$ with $\|x_n\|\to 0$ and $\delta(x_n) \xrightarrow{m} \eta$.  Then
\begin{align*}
0 = \delta([x_n, h])&= (x_n \otimes 1^\op)\delta(h)+(1\otimes h^\op)\delta(x_n) - (1\otimes x_n^\op)\delta(h)-(h\otimes 1^\op)\delta(x_n)\\
&\xrightarrow{m} (1\otimes h^\op-h\otimes 1 )\eta, 
\end{align*}
and since $h$ is diffuse $(1\otimes h^\op-h\otimes 1 )$ is not a zero-divisor in $\ms U$;  hence $\eta=0$. We now claim that $\delta$ has to be inner on $Z(M)=L\ZZ$. If this were not the case, then, by Lemma \ref{lemma:vanishing-change-of-coeff}, there exists a non-inner derivation $\delta' \colon L\ZZ\to \ms{U}(L\ZZ\htens L\ZZ)$. By what was just proven, $\delta'$ is norm-measure continuous and therefore, by Theorem \ref{sot-mt-cont},  also continuous for the strong operator topology.  Because of this, the restriction of $\delta'$ to the complex group algebra $\CC[\ZZ]$  has to be non-inner, contradicting the fact that $0=\beta_1^{(2)}(\ZZ)=\beta_1^{(2)}(\CC[\ZZ],\tau)$ (cf.~\cite[Proposition 2.3]{CS} and \cite[Theorem 0.2]{cheeger-gromov}).  Hence there exists $\xi\in \ms U$ such that $\delta$ agree with $\delta_\xi:=[\cdot, \xi]$ on $Z(M)$. The difference $\delta-\delta_\xi$ therefore vanishes on $Z(M)$ and by Lemma \ref{lemma:vanishing-on-the-normalizer} it has to vanish on every unitary in $M$. Since 
the unitaries span $M$ linearly, we conclude that $\delta$ is globally inner.
\end{proof}
\begin{Rem}
By combining Lemma \ref{lemma:vanishing-on-the-normalizer} and Proposition \ref{beta-one-of-diffuse} with the automatic strong continuity from Theorem \ref{sot-mt-cont}, one may at this point easily deduce the conclusion of \cite[Theorem 6.4]{thom2008l2}; namely that the first continuous $L^2$-Betti number vanishes for non-prime von Neumann algebras as well as for von Neumann algebras admitting a diffuse Cartan subalgebra. Since this will also follow from the more general vanishing results obtained in Section \ref{compression-subsection}, we shall not elaborate further at this point.
\end{Rem}

\subsection{Factors with property (T)}\label{prop-T-section}
If $\Gamma$ is a countable discrete group with property $\T$ it is well known that its first $L^2$-Betti number vanishes. This observation dates back to the work of Gromov \cite{gromov-asymptotic-invariants}, but the first complete  proof was given by Bekka and Valette in \cite{bekka-valette-group-cohomology}. Applying the recent techniques from\cite[Theorem 2.2]{peterson-thom}, this can now be deduced easily from the Delorme-Guichardet theorem (see e.g.~\cite{kazhdan's-property-t}), which characterizes property (T) of $\Gamma$ in terms of vanishing of its first cohomology groups. 
The notion of property $\T$ for II$_1$-factors was introduced  by Connes and Jones in \cite{connes-jones} and in \cite{peterson-T} Peterson proved a version of the Delorme-Guichardet theorem in this context:

\begin{Thm}[{\cite[Theorem 0.1]{peterson-T}}]\label{thm:jesse's}
 Let $M$ be a finite factor with separable predual. Then the following conditions are equivalent:
\begin{enumerate}
 \item $M$ has property $\T$;
 \item there exists a weakly dense $\ast$-subalgebra $M_0\subset M$ which is countably generated as a vector space and such that every densely defined $L^2$-closable derivation from $M$ into a Hilbert $M$-$M$-bimodule whose domain contains $M_0$ is inner.
\end{enumerate}
\end{Thm}

In this section we apply Peterson's result to prove the following von Neumann algebraic version of the classical group theoretic result mentioned above.
\begin{Thm}\label{beta-1-of-T-factor}
Let $M$ be a $\twoone$-factor with separable predual. If $M$ has property $\T$ then $\eta_1^{(2)}(M,\tau)=0$.
\end{Thm}
\begin{proof}
Since $M$ has property $\T$, we obtain from Theorem \ref{thm:jesse's} a dense $*$-sub\-alge\-bra $M_0\subset M$ such that $M_0$ is countably generated as a vector space and such that any derivation from $M_0$ into a Hilbert $M$-bimodule $H$, which is closable as an unbounded operator from $L^2(M)$ to $H$, is inner.  We first observe that
\begin{align*}
\dim_{M\htens M^\op} \Der_c(M,\ms U)=\dim_{M\htens M^\op}\{\delta\in \Der_c(M,\ms U) \mid \delta(M_0)\subset M\htens M^\op\}.
\end{align*} 
To see this, it suffices by Sauer's local criterion \cite[Theorem 2.4]{sauer-betti-for-groupoids} to prove that for each continuous derivation $\delta\colon M \to \ms U$ and each $\eps>0$ there exists a projection $p\in M\htens M^\op$ such that $\tau(p^\perp)\leqslant\eps$ and $\delta(-)p$ maps $M_0$ into $M\htens M^\op$. Choose a countable linear basis $(e_n)_{n=1}^\infty$ for $M_0$. Since each $\delta(e_n)$ is affiliated with $M\htens M^\op$ there exists a projection $p_n\in M\htens M^\op$ such that $\tau(p_n^\perp)\leqslant\frac{\eps}{2^n}$ and such that $\delta(e_n)p_n\in M\htens M^{\op}$. The projection $p:=\bigwedge_n p_n$ therefore satisfies the requirements. 
We now have to prove that
\[
\dim_{M\htens M^\op}\{\delta\in \Der_c(M,\ms U) \mid \delta(M_0)\subset M\htens M^\op\}\leqslant 1,
\]
and we will do so by proving that a continuous derivation $\delta\colon M\to \ms U$ for which $\delta(M_0)\subset M\htens M^\op$ has to be inner. We claim that it suffices to prove that $\delta$ is $L^2$-closable from $M_0$ to $L^2(M\htens M^\op)$. Indeed, if this is the case, then by Peterson's result there exists a vector $\xi\in L^2(M\htens M^\op)$ such that 
\[
\delta(a)=(a\otimes 1)\xi - (1\otimes a^\op)\xi \ \text{ for all } a \in M_0.
\]
Considering $\xi$ as an operator in $\ms U$, we get that it implements $\delta$ on $M_0$ and hence by Theorem \ref{sot-mt-cont} it implements $\delta$ on all of $M$.
Thus, our task is  to show that $\delta\colon M_0 \to L^2(M\htens M^\op)$ is $L^2$-closable. Let therefore $x_n\in M_0$ and assume that $x_n \xrightarrow{2} 0$ and $\delta(x_n)\xrightarrow{2} \eta$. Since convergence in $2$-norm implies convergence in measure (Lemma \ref{Lp-lem}) and since $\delta$ is continuous from the measure topology on $M$ to the measure topology on $\ms U$ (Theorem \ref{sot-mt-cont}) we obtain that $\delta(x_n)\overset{m}{\to} 0$ as well as $\delta(x_n)\overset{m}{\to}\eta$; hence $\eta=0$ as the measure topology is Hausdorff.
\end{proof}

\subsection{Further remarks}
In this section we collect a few observations concerning the first continuous $L^2$-Betti number. Some of them are easily derived from the results in \cite{thom2008l2}, but since they are also direct consequences of Theorem \ref{sot-mt-cont} we include them here for the sake of completeness.

\begin{Prop}\label{weakly-dense-prop}
Let $A\subset M$ be a weakly dense $*$-subalgebra. Then 
 \[
 \eta_1^{(2)}(M,\tau)\leqslant \eta_1^{(2)}(A,\tau)\leqslant \beta_1^{(2)}(A,\tau).
 \]
and if $A$ is a \Cs-algebra we have  $\eta_1^{(2)}(M,\tau)= \eta_1^{(2)}(A,\tau)$.
\end{Prop}
Note that the inequalities in Proposition \ref{weakly-dense-prop} are contained in \cite[Theorem 6.2]{thom2008l2}, but that the equality  when $A$ is a \Cs-algebra does not directly follow from this. Compare also with \cite[Theorem 4.6]{thom2008l2}.

\begin{proof}
By Theorem \ref{sot-mt-cont}  the map $H_c^1(M,\ms U) \longrightarrow H_c^1(A, \ms U)$ induced by restriction is injective in general and an isomorphism when $A$ is a \Cs-algebra. That $\eta_1^{(2)}(A,\tau)\leqslant \beta_1^{(2)}(A,\tau)$  is clear, as we have an inclusion $\Der_c(A,\ms{U}) \subset \Der(A,\ms U)$ for any algebra $A$. 
\end{proof}
To illustrate the usefulness of the above result, we record the following  consequences.

\begin{itemize}
\item[1)] If $\Gamma$ is a discrete countable group with $\beta_1^{(2)}(\Gamma)=0$ then also $\eta_1^{(2)}(L\Gamma,\tau)=0$. By \cite[Proposition 2.3]{CS}, $\beta_1^{(2)}(\Gamma)=\beta_1^{(2)}(\CC\Gamma, \tau)$ and the claim therefore follows from Proposition \ref{weakly-dense-prop}. In particular, the first continuous $L^2$-Betti number of the hyperfinite factor $R$ vanishes since $R\simeq L\Gamma$ for any amenable icc group $\Gamma$.
The result about the hyperfinite factor can also be obtained directly from the definition of hyperfiniteness by realizing $R$ as the von Neumann algebraic direct limit of matrix algebras, for which it is also known \cite[Example 6.9]{dim-flatness} that the (non-continuous) $L^2$-Betti numbers  of the corresponding \emph{algebraic} direct limit vanishes.

\item[2)]  It is well known \cite{bekka-valette-group-cohomology} that when $\Gamma$ is a discrete, countable group with property $\T$ then $\bet_1(\Gamma)=0$ and hence also $\eta_1^{(2)}(L\Gamma,\tau)=0$. Thus, in the case of factors arising from discrete groups, Theorem \ref{beta-1-of-T-factor} can be deduced immediately.

\item[3)] For the von Neumann algebra $L^\infty(O_n^+)$ associated with the free orthogonal quantum group $O_n^+$ we have $\eta_1^{(2)}(L^\infty(O_n^+),\tau)=0$. 
Denoting by $\Pol(O_n^+)$ the canonical dense Hopf $*$-algebra in $L^\infty(O_n^+)$, it is known  that  $\beta_1^{(2)}(\Pol(O_n^+),\tau)=0$ (see \cite{coamenable-betti} for the case $n=2$ and \cite{vergnioux-paths-in-cayley} for the case $n\geqslant3$) and hence, by Proposition \ref{weakly-dense-prop}, $\eta_1^{(2)}(L^\infty(O_n^+),\tau)=0$.

\end{itemize}

Since $\eta_1^{(2)}(-)$ measures the dimension the space of continuous derivations it follows from the results already proven that this number is finite for von Neumann algebras that are finitely generated. Since we will use this repeatedly in the sequel, where a concrete upper bound will be of importance, we single this out by means of the following lemma.
\begin{Lemma}\label{fg-lem}
If $M$ is generated as a von Neumann algebra by $n$ selfadjoint elements then $\eta_1^{(2)}(M,\tau)\leqslant n-1$.
\end{Lemma}
\begin{proof}

If $M$ is generated by $n$ selfadjoint elements $x_1,\dots, x_n$ then  the complex subalgebra $A$ generated by $\{1, x_1,\dots, x_n \}$ is a dense unital $*$-subalgebra in $M$ and by Proposition \ref{weakly-dense-prop} any continuous derivation $\delta\colon M \to \ms U$ is uniquely determined by its values on $A$. From the derivation property it follows that $\delta$ is already completely determined on its values on the generators $x_1,\dots, x_n$ and hence we get
$\eta_1^{(2)}(M,\tau)\leqslant n-1$ as desired. 
\end{proof}

\subsection{The compression formula in continuous cohomology}\label{compression-subsection}

Recall from \cite[Theorem 2.4]{CS} that the algebraic $L^2$-Betti numbers scale quadratically when passing to corners; more precisely if $M$ is a finite factor and $p\in M$ is a non-zero projection then $\beta_n^{(2)}(pMp, \tau_p)=\tau(p)^{-2}\beta^{(2)}_n(M,\tau)$, where $\tau_p$ denotes the restriction of $\tau$ to the corner $pMp$ rescaled with $\tau(p)^{-1}$. In this section we prove that the same holds true for the first continuous $L^2$-Betti number 
and, as a byproduct, provide a cohomological proof of the scaling formula for the first algebraic $L^2$-Betti number as well. 
\begin{Thm}\label{scaling-for-eta}
Let $M$ be a $\twoone$-factor with trace-state $\tau$ and let $p\in M$ be a non-zero projection. Then $\eta_1^{(2)}(pMp, \tau_p)=\frac{1}{\tau(p)^2}\eta_1^{(2)}(M,\tau)$.
\end{Thm}\label{thm:scaling-formula}
\begin{proof}
Denote $p\tens p^\op\in M\htens M^\op$ by $q$ and consider the right $qM\htens M^\op q$-linear map
\begin{align*}
\Phi_p\colon \Der(M,\ms U q) &\longrightarrow \Der(pMp, q\ms U q)\\
\delta & \longmapsto (p\otimes p^\op) \cdot \delta|_{pMp}
\end{align*}
Note that this map induces  a map $\Phi_p\colon H^1(M,\ms Uq) \to H^1(pMp,q\ms Uq)$ on the ordinary cohomology as well as a map on the continuous cohomology $\Phi_p\colon H^1_c(M,\ms U q) \to H^1_c(pMp,\ms U q)$, since continuity of a derivation is preserved by construction, and an inner derivation implemented by $\xi\in \ms U q$ is mapped to the inner derivation on $pMp$ implemented by $q\xi$. We also note that the restriction maps $\Phi_*$ are compatible with the order structure: if $r,s\in \Proj(M)$ and $r\leqslant s$ then $\Phi_r= \Phi_r\circ\Phi_s$. Our aim is to prove that the map $\Phi_p$ is an isomorphism of right $qM\htens M^\op q$-modules on the level of continuous 1-cohomology. Once this is established, the result follows from the general cut-down formula for the dimension function (see e.g.~\cite[Lemma A.15]{KPV12}) since
\begin{align*}
\eta_{1}^{(2)}(pMp,\tau_p)&:=\dim_{qM\htens M^\op q} H^1_c\left(pMp, \ms U(pMp\htens (pMp)^\op)\right)\\
&= \dim_{qM\htens M^\op q} H^1_c(pMp, q\ms Uq) \\
&=\dim_{qM\htens M^\op q} H^1_c(M,\ms Uq) \\
&=\dim_{qM\htens M^\op q} H^1_c(M,\ms U)q \\
& = \frac{1}{(\tau\tens\tau^\op)(q)}\dim_{M\htens M^\op}H^1_c(M,\ms U)\\
&=\frac{1}{ \tau(p)^2} \eta_1^{(2)}(M,\tau). 
\end{align*}
By construction, $\Phi_p$  is right $q(M\htens M^\op)q$-linear so we only have to provide the inverse to $\Phi_p \colon H^1(M,\ms Uq) \to H^1(pMp,q\ms Uq) $ and show that it maps $ H_c^1(pMp,q\ms Uq)$ to $H_c^1(M,\ms Uq)$. To this end, choose $n$ to be the smallest integer such that $n\tau(p)\geqslant 1$ and choose orthogonal projections $p_1,\dots, p_{n}\in M$ summing to $1_M$ such that $p_1,p_2,\dots, p_{n-1}$ 
are equivalent to $p$ and $p_n$ is equivalent to a subprojection $f$ of $p$. We furthermore may, and will, assume that $p_1=p$.
This choice provides us with a $*$-isomorphism $s\mb M_n(pMp)s \cong M$ where $s\in \MM_n(pMp)=\MM_n(\CC)\tens pMp$ is the projection $\sum_{i=1}^{n-1} v_{ii}\tens p +v_{nn}\tens f$. Here, and in what follows, we denote by $\{v_{ij}\}_{i,j=1}^n$ the standard matrix units in $\MM_n(\CC)$. In the sequel we will suppress this isomorphism and simply identify $M$ and $s\MM_n(pMp)s$.\\

Now we define two maps. The first one is ``induction-to-matrices'':
$$
\ind_n\colon \Der(pMp,\ms U q) \to \Der\left(\mb M_n(pMp), \mb M_n(\mb C)\otimes \mb M_n(\mb C)^\op\otimes q\ms U q\right)
$$
given by
$$
\ind_n(\delta)(x) = \sum_{i,j=1}^n (v_{i,1}\otimes v_{1,j}^\op)\tens  \delta((v_{1,i}\tens p) x (v_{j,1}\tens p)).
$$
A direct computation verifies  that $\ind_n(\delta)$ is indeed a derivation when $\mb M_n(\mb C)\otimes \mb M_n(\mb C)^\op\otimes q\ms U q$ is endowed with the natural $\MM_n(\CC)\tens pMp$-bimodule structure given by
\[
(a\tens x)\tens (b\tens y)^\op.T:= \left( a\tens b^\op \tens (x\tens y^\op)\right)T, \ a,b\in \MM_n(\CC), x,y\in pMp, T\in q\ms U q.
\]
The map $\ind_n(\delta)$ descends to both cohomology and continuous cohomology since an inner derivation implemented by $\xi \in q\ms U q$ maps to the inner derivation implemented by $\sum_{i=k}^n v_{k1}\tens v_{1k}^\op\tens \xi$
and since $\ind_n(-)$ clearly maps continuous derivations to continuous derivations.\\

The second map is the compression map with respect to the projection $s$:
\begin{align*}
\Phi_s\colon \Der\big(\mb M_n(pMp), \mb M_n(\mb C)\otimes \mb M_n(\mb C)^\op\otimes  q\ms Uq\big)& \longrightarrow \Der(M,\ms U q)\\
\delta & \longmapsto s\otimes s^\op \cdot \delta|_{s\mb M_n(pMp)s}
\end{align*}
Note that this restriction map indeed maps to $\Der(M,\ms U q)$: in the matrix picture $p$ identifies with $v_{11}\tens p$ which is a subprojection of $s$ and hence
\begin{align*}
\ms U q &= \ms U \Big ( s\tens s^\op \big(\MM_n(\CC)\tens pMp\big)\htens \big(\MM_n(\CC)^\op\tens (pMp)^\op\big) s\tens s^\op  \Big)q\\
&=(s\tens s^\op) \ms U \Big( (\MM_n(\CC)\tens pMp)\htens (\MM_n(\CC)^\op\tens (pMp)^\op)    \Big)\\
&= (s\tens s^\op).  (\MM_n(\CC)\tens \MM_n(\CC)^\op \tens q\ms U q).
\end{align*}
We now claim that $\Phi_s \circ \ind_n $ is the inverse of $\Phi_p$ on the level of (continuous) cohomology. One composition can be easily computed using the order compatibility of the restriction maps:
$
\Phi_p\circ \Phi_s\circ \ind_n = \Phi_p\circ \ind_n
$
which is the identity map even at the level of derivations. To see this, consider $x\in pMp$ and recall that in the matrix picture $p\in M$ identifies with the projection $v_{11}\tens p$. Thus,
\begin{align*}
\Phi_p\circ \ind_n(\delta) (x)&=  (v_{11}\tens p)\tens (v_{11}\tens p)^\op\left(\sum_{i,j=1}^n v_{i1}\tens v_{1j}^\op \tens \delta\big((v_{i1}\tens p)x(v_{j1}\tens p)  \big)\right)\\
&=v_{11}^{\phantom{\op}}\tens v_{11}^\op \tens \delta(x) =\delta(x).
\end{align*}
Next we have to compute $\Phi_s\circ \ind_n\circ \Phi_p$. We start with a derivation $\delta \in \Der(M,\ms U q)$. Consider the the following two systems of matrix units $\{v_{ij}\tens (p-f)\}_{i,j=1}^{n-1}$ and $\{v_{ij}\tens f\}_{i,j=1}^n$ in $M$ and denote by $A_0$ the $*$-algebra they generate. This is a finite dimensional\footnote{The map $\mb M_{n-1}(\mb C)\oplus \mb M_{n}(\mb C)\ni (a,b)\mapsto a\otimes (p-f) + b\otimes f$ is surjective and an isomorphism in the generic case when $f\neq 0$ and $f\neq p$.}  $C^*$-algebra and hence  $\beta_1^{(2)}(A_0,\tau|_{A_0})=0$ by \cite[Proposition 2.9]{CS}. By Lemma \ref{lemma:vanishing-change-of-coeff}, the restriction of $\delta$ to $A_0$ is therefore inner, so by subtracting an inner derivation we may assume that $\delta$ vanishes on $A_0$.  Hence for all $a,b\in A_0$ and $x\in M$ we have
$
\delta(axb)=(a\tens b^\op)\delta(x).
$
In particular,  $\Phi_p(\delta)=\delta|_{pMp}$ since $p=v_{11}\tens p \in A_0$. Thus, splitting the unit in $M$ as $\sum_{i=1}^{n-1} v_{ii} \tens p +v_{nn}\tens f$ we have 
\begin{align*}
\delta(x)&= \sum_{i,j=1}^{n-1} \delta\big((v_{ii}\tens p)x(v_{jj}\tens p)\big) +\sum_{i=1}^{n-1}\delta\big((v_{ii}\tens p)x(v_{nn}\tens f)\big) +\\
&+\sum_{i=1}^{n-1}\delta \big((v_{nn}\tens f)x(v_{ii}\tens p)\big) + \delta\big((v_{nn}\tens f) x(v_{nn}\tens f)\big)\\
    &=\sum_{i,j=1}^{n-1} \delta\big((v_{i1}\tens p)(v_{1i}\tens p)x(v_{j1}\tens p) (v_{1j}\tens p)\big) +\\
    &+\sum_{i=1}^{n-1}  \delta\big((v_{i1}\tens p)(v_{1i}\tens p)x(v_{n1}\tens p)(v_{1n}\tens f)\big) +\\
    &+\sum_{i=1}^{n-1}  \delta \big( (v_{n1}\tens f)(v_{1n}\tens p)x(v_{i1}\tens p)(v_{1i}\tens p)\big)+\\ 
    &+ \delta\big((v_{n1}\tens f)(v_{1n}\tens p) x(v_{n1}\tens p)(v_{1n}\tens f)\big)\\
    &=\sum_{i,j=1}^{n-1} (v_{i1}\tens p)\tens(v_{1j}\tens p)^\op \delta \big((v_{1i}\tens p)x(v_{j1}\tens p)\big) +\\
    &+\sum_{i=1}^{n-1} (v_{i1}\tens p)\tens (v_{1n}\tens f)  \delta\big((v_{1i}\tens p)x(v_{n1}\tens p)\big) +\\
    &+\sum_{i=1}^{n-1}(v_{n1}\tens f)\tens(v_{1i}\tens p)^\op  \delta\big((v_{1n}\tens f)x(v_{i1}\tens p)\big)+\\ 
    &+ (v_{n1}\tens f)\tens (v_{1n}\tens f)\delta\big((v_{1n}\tens p) x(v_{n1}\tens p)\big)\\
&= (s\tens s^\op). \left( \sum_{i,j=1}^n (v_{i1}\tens v_{1j}^\op) \tens \delta((v_{1i}\tens p)x(v_{j1}\tens p))\right)\\
&=(s\tens s^\op).\ind_{n}(\delta|_{pMp}) (x)\\
&=(s\tens s^\op). \ind_n(\Phi_p(\delta))(x)\\
&= \Phi_s\circ \ind_n \circ \Phi_p(\delta) (x),
    \end{align*}
as desired.    
\end{proof}
\begin{Rem}\label{scaling-rem}
More generally, for any $t>0$ the $t$-th amplification of $M$ is defined as $r\MM_n(M)r$ where $n= \lfloor t \rfloor +1$ and $r\in \MM_n(M)$ is a projection with $(\tr_n\tens \tau)(r)=t/n$. We note that the scaling formula holds true in this generality since applying it first to $M$ considered as a corner in $\MM_n(M)$ yields $\eta_1^{(2)}(M)=n^{-2}\eta_1^{(2)}(\MM_n(M))$ and applying it once more with respect to the projection $r$ therefore gives
\[
\eta_1^{(2)} (M_t)=(\tr_n\tens \tau)(r)^{-2}\eta_1^{(2)}(\MM_n(M))= \frac{n^2}{t^2}n^2\eta_1^{(2)}(M)=t^{-2}\eta_1^{(2)}(M).
\]

\end{Rem}

The isomorphism provided in the proof of Theorem \ref{thm:scaling-formula} is clearly also an isomorphism on the algebraic level so along the way we also proved the following special case of \cite[Theorem 2.4]{CS}.
\begin{Por}
If $M$ is a $\twoone$-factor and $t>0$ then $\beta_1^{(2)}(M_t)=t^{-2}\beta_1^{(2)}(M)$.
\end{Por}

\begin{Cor}\label{fundamental-grp-cor}
If $M$ is a finitely generated $\twoone$-factor with  non-trivial fundamental group then $\eta_1^{(2)}(M,\tau)=0$.
\end{Cor}

\begin{proof}
If $M$ is generated by $n$ elements $x_1,\dots, x_n$ then by extracting the real and imaginary part of these generators we get $2n$ selfadjoint generators $a_1,\dots, a_n, b_1,$ $\dots, b_n$ and from Lemma \ref{fg-lem} it follows that $\eta_1^{(2)}(M)\leqslant 2n<\infty$. Picking a non-trivial projection $p\in M$ such that $M\cong pMp$\footnote{Since $M$ is a factor this isomorphism must intertwine the trace $\tau$ with $\tau_p$.}
we conclude from the scaling formula that
\[
\eta_1^{(2)}(M,\tau)=\eta_1^{(2)}(pMp,\tau_p)=\frac{1}{\tau(p)^2}\eta_1^{(2)}(M,\tau),
\]
and since we just argued that $\eta_1^{(2)}(M,\tau)<\infty$ this forces $\eta_1^{(2)}(M,\tau)=0$.
\end{proof}
Recently, Shen \cite{She05} introduced the \emph{generator invariant} $\mathcal{G}(M)$ and proved that $\mathcal{G}(M)<\frac14$ implies that $M$ is singly generated. A further study of the generator invariant, as well as its hermitian analogue $\mathcal{G}_\sa(M)$, was undertaken in \cite{DSSW08} where the authors, inter alia, prove a scaling formula for the invariant under the passage to corner algebras.  This scaling formula implies that the class of $\twoone$ factors with vanishing generator invariant\footnote{At the time of writing, no example of a $\twoone$ factor with non-vanishing generator invariant is known.} is stable under passing to corners.
As a consequence, we obtain the following result showing, yet again, that the first continuous $L^2$-Betti number vanishes on a large class of $\twoone$ factors.

\begin{Cor}\label{bounded-by-G-cor}
For any $\twoone$ factor $M$ we have $\eta_1^{(2)}(M)\leqslant  \mathcal{G}_\sa(M)$. In particular, factors that are either non-prime, admits a Cartan or has property Gamma has vanishing first continuous $L^2$-Betti number.
\end{Cor}
The inequality in  Corollary \ref{bounded-by-G-cor} can be deduced from the more general result \cite[Corollary 5.12]{DSSW08}, but since the result there is stated without proof we find it worthwhile to include the short argument below. Note also that the vanishing  of $\eta_1^{(2)}(-)$ in the non-prime and Cartan case was already proved by Thom in \cite{thom2008l2} and that in the special case of a \emph{group} von Neumann algebra,  the result about property Gamma factors can be deduced from \cite[Theorem 1.2]{peterson-L2-rigidity} and the general inequality $\eta_1^{(2)}(L\Gamma)\leqslant \beta_1^{(2)}(\Gamma)$.

\begin{proof}
We first prove that if $k$ is any integer then if $\mathcal{G_\sa}(M)< k$ then $\eta_1^{(2)}(M)\leqslant k$. To see this, just note that by \cite[Theorems 3.1 \& 5.5]{DSSW08}
we have that $M$ is generated by $k+1$ selfadjoint elements and by Lemma \ref{fg-lem} this implies that $\eta_1^{(2)}(M)\leqslant  k$. The inequality is trivial when $\mc G_{\sa}(M)=\infty$, so assume that $\mc G_\sa(M)<\infty$ and put $t_n:=\sqrt{\mc G_\sa(M)+\frac1n}$. Then by the scaling formula for $\mc_\sa$ \cite[Corollary 5.6]{DSSW08} we have $\mc G_\sa(M_{t_n})<1$ and hence by what was just proven also $\eta_1^{(2)}(M_{t_n})\leqslant 1$. By Theorem \ref{scaling-for-eta} (see also Remark \ref{scaling-rem}) we therefor have 
\[
\eta_1^{2}(M)=t_n^{2}\eta(M_{t_n})\leqslant t_n^2\underset{n\to \infty}{\longrightarrow} \mc G_\sa(M).
\]
That $\eta_1^{(2)}(M)=0$ when $M$ is non-prime, has a Cartan subalgebra or has property Gamma follows from the formula $\mc G_sa(M)=\frac12 \mc G_\sa(M)$ \cite[Theorem 5.5]{DSSW08} in conjunction with \cite[Section 6]{She05} where it is shown that $\mathcal{G}(M)=0$ under the aforementioned hypotheses.
\end{proof}

\begin{Cor}\label{fundamenal-grp-cor}
The first continuous $L^2$-Betti number vanishes on any class of singly generated $\twoone$-factors that is stable under passing to corners.

\end{Cor}

\begin{proof}
Let $\mathcal C$ be such a class of $\twoone$ factors and note that $\eta_1^{(2)}(-)$ is bounded by 1 on $\mc C$. Let $M\in \mc C$ be given and choose a sequence of non-trivial projections $p_n\in M$ with $\tau(p_n)\to 0$. Since $\mc C$ is stable under passing to corners we obtain 
\[
1\geqslant \eta_1^{(2)}(p_nMp_n,\tau_{p_n})=\frac{1}{\tau(p_n)^2}\eta_1^{(2)}(M,\tau). 
\]
and hence
\[
\eta_1^{(2)}(M,\tau)\leqslant  \tau(p_n)^2\underset{n\to \infty}{\longrightarrow} 0.
\]
\end{proof}
Note that above corollary has the following curious consequence: If the notorious generator problem has a positive solution (i.e.~every $\twoone$ factor is singly generated) then the first  continuous $L^2$-Betti number vanishes globally on the class of $\twoone$ factors. \\

We end this section with a  result regarding interpolated free group factors which shows that the first continuous $L^2$-betti number is ``linear in the number of generators''.
This is another consequence of the scaling formula and the proof is verbatim the same as the corresponding proof regarding Shen's generator invariant given in \cite{DSSW08}. We include it below for the sake of completeness.

\begin{Prop}
There exists $a\in[0,1]$ such that $\eta_1^{(2)}(L\FF_r)=a(1-r)$ for every $r\in]1,\infty[$. In particular, if $\eta_1^{(2)}(L\FF_2)>0$ then the interpolated free group factors $L\FF_r$ are pairwise non-isomorphic for $r\in ]1,\infty[$.
\end{Prop}
\begin{proof}
Denote by $f\colon ]0,\infty[\to \mathbb R$ the function $r\mapsto \eta_1^{(2)}(L\FF_{r+1})$ and recall \cite{dykema-interpolated} that for $r>1$ and $\lambda>0$ the interpolated free group factors satisfy the scaling formula
$
L(\FF_{1+\frac{r-1}{\gamma^2}})=L(\FF_r)_\gamma.
$
Hence
\[
f\left(\frac{r-1}{\lambda^2}\right)=\eta_1^{(2)}\left(\FF_{1+\frac{r-1}{\gamma^2}}\right)=\frac{1}{\gamma^2}\eta_1^{(2)}\left(L\FF_r\right)= \frac{1}{\gamma^2}f(r-1).
\]
The map $f$ therefore satisfies $f(sr)=sf(r)$ for all $s,r>0$ and thus $f(r)=rf(1)=r\eta_1^{(2)}(L\FF_2)$; hence $a:=\eta_1^{(2)}(L\FF_2)$ does the job and since $\eta_1^{(2)}(L\FF_2)\leqslant  \beta_1^{(2)}(\FF_2)=1$ we have $a\in [0,1]$. The final statement concerning non-isomorphism follows trivially from this.

\end{proof}

\bibliographystyle{alpha}

\end{document}